\documentclass[12pt]{amsart}

\usepackage{amssymb}
\usepackage{color}

\usepackage{enumitem}

\makeatletter
\@namedef{subjclassname@2020}{%
  \textup{2020} Mathematics Subject Classification}
\makeatother

\newtheorem{theorem}{Theorem}[section]
\newtheorem{corollary}[theorem]{Corollary}
\newtheorem{lemma}[theorem]{Lemma}
\newtheorem{proposition}[theorem]{Proposition}

\theoremstyle{definition}
\newtheorem{definition}[theorem]{Definition}

\numberwithin{equation}{section}

\begin{document}


\baselineskip=17pt

\title[Mapped coercivity for the stationary Navier-Stokes equations]{Mapped coercivity for the stationary Navier-Stokes equations and their finite element approximations}

\author{Roland Becker}
\address[R. Becker]{Department of Mathematics, Universit\'e de Pau et de l'Adour (UPPA), Avenue de l'Université, BP 1155, 64013 Pau CEDEX, France}
\email{becker@univ-pau.fr}

\author{Malte Braack}
\address[M. Braack]{Mathematical Seminar, Kiel University, Heinrich-Hecht-Platz 6, 24098 Kiel}
\email{braack@math.uni-kiel.de}

\date{\today}

\begin{abstract}
This paper addresses the challenge of proving the existence of solutions for nonlinear equations in Banach spaces, 
focusing on the Navier-Stokes equations and discretizations of thom. 
Traditional methods, such as monotonicity-based approaches and fixed-point theorems, often face limitations in handling general nonlinear operators 
or finite element discretizations. 
A novel concept, mapped coercivity, provides a unifying framework to analyze nonlinear operators through a continuous mapping. 
We apply these ideas to saddle-point problems in Banach spaces, emphasizing both infinite-dimensional formulations and finite element discretizations. 
Our analysis includes stabilization techniques to restore coercivity in finite-dimensional settings, ensuring stability and existence  of solutions. 
For linear problems, we explore the relationship between the inf-sup condition and mapped coercivity, using the Stokes equation as a case study. 
For nonlinear saddle-point systems, we extend the framework to mapped coercivity via surjective mappings,
enabling concise proofs of existence  of solutions for various stabilized Navier-Stokes finite element methods. 
These include Brezzi-Pitkäranta, a simple variant, and local projection stabilization (LPS) techniques, with extensions to convection-dominant flows.
The proposed methodology offers a robust tool for analyzing nonlinear PDEs and their discretizations, 
bypassing traditional decompositions and providing a foundation for future developments in computational fluid dynamics.

\end{abstract}

\subjclass[2020]{Primary ; Secondary }

\maketitle

%


\def\beqa{\begin{eqnarray*}}
\def\eeqa{\end{eqnarray*}}
\def\beqal{\begin{eqnarray}}
\def\eeqal{\end{eqnarray}}
\def\bu{\boldsymbol u}
\def\bv{\boldsymbol v}
\def\bob{\boldsymbol f}
\def\bw{\boldsymbol w}
\def\bV{\boldsymbol V}
\def\bD{\boldsymbol D}
\def\bX{\boldsymbol X}
\def\bphi{\boldsymbol\phi}
\def\nR{\mathbb R}
\def\div{{\rm div}\,}
\def\norm#1{\|#1\|_{L^2(\Omega)}}
\def\normX#1{\|#1\|_{X}}
\def\normXS#1{\|#1\|_{X'}}
\def\normT#1{\|#1\|_{L^\infty(T)}}
\def\normSY#1{\|#1\|_{Y'}}
\def\normY#1{\|#1\|_{Y}}
\def\normV#1{\|#1\|_{V}}
\def\normVS#1{\|#1\|_{V'}}
\def\normQ#1{\|#1\|_{Q}}
\def\normQS#1{\|#1\|_{Q'}}
\def\normW#1{\|#1\|_{Y',Y}}
\def\normY#1{\|#1\|_{V'\leftarrow V}}
\def\vector#1{\left(\begin{array}{c}#1\end{array}\right)}
\def\l2#1#2{(#1,#2)_{L^2(\Omega)}}
\def\endproof{\hfill$\square.$}
\renewcommand{\epsilon}{\varepsilon}
\newcommand{\Set}[1]{ \left\{#1\right\}}
\newcommand{\SetDef}[2]{\left\{#1\;\middle|\;#2\right\}} 
\newcommand{\Rest}[2]{{#1}_{|_{#2}}}
\newcommand{\transpose}[1]{#1^{\mathsf{T}}} 
\newcommand{\transposeInv}[1]{#1^{\mathsf{-T}}} 
\newcommand{\scp}[2]{\langle#1,#2\rangle}
\newcommand{\Norm}[1]{\left\|#1\right\|}
\definecolor{myred}{rgb}{0.58, 0.066, 0.}
\newcommand{\red}[1]{\textcolor{myred}{#1}}

\section{Introduction}

The available methods for proving the existence of solutions to nonlinear equations in Banach spaces are fairly limited. 
For example, the foundational theory by Browder and Minty hinges on monotonicity, 
a condition that applies only to specific types of equations. 
Meanwhile, fixed-point theorems such as those by Leray and Schauder require reformulating 
the equation into a fixed-point form that meets the requirements of suitable function spaces, 
often involving carefully chosen preconditioners.

In the classical approach to proving the existence of solutions for the Navier-Stokes equations, 
the system’s saddle-point structure is exploited. This allows for a decomposition into two sub-problems: 
a nonlinear equation for the velocity in the divergence-free function space, controlled by a coercive operator, 
and a linear equation for the pressure. However, many finite element discretizations complicate this two-step approach,
as the discrete velocity field may fail to be divergence-free. Consequently, a fully coupled system must be addressed,
as demonstrated by Tobiska and Verfürth in their work on the SUPG discretization \cite{TobiskaVerfuerth1996}. 
Their proof confirms the existence of discrete solutions for right-hand sides with bounded norms, 
where the bound is determined by the mesh size and spatial dimension.

A recent advance in solution theory, introduced in \cite{BeckerBraack2024}, proposes the concept of 
{\em mapped coercivity}. 
This approach aims to streamline and enhance the analysis of nonlinear operators in Banach spaces. 
An essential aspect of this theory is the use of a continuous mapping $\Phi:X\to X$
to establish a generalized form of coercivity. In specific cases, this map can be linear, giving rise to {\em linearly mapped coercivity}. This paper seeks to apply this new concept within the context of saddle-point problems in a Banach space of the form $X=V\times Q$,
with a focus on the Navier-Stokes equations in both their infinite-dimensional formulation and finite element discretizations.

We begin by examining linear saddle-point systems in the context of {\em linearly mapped coercivity}. 
When this property does not naturally extend to finite-dimensional formulations, 
stabilization terms may be introduced to restore coercivity, 
thereby ensuring stability of the finite-dimensional approximations. 
This approach is broadly applicable, as it does not rely on the structure of an inner product. 
As a result, the general theory of saddle-point problems remains valid in Banach spaces. 
The proposed methodology offers a versatile framework for proving the existence of solutions 
to a wide range of nonlinear partial differential equations and their discretizations.

To explore concrete partial differential equations, we analyze a generalized Stokes system within a 
Hilbert space setting and review several established stabilization techniques. 
We demonstrate that the corresponding operators exhibit {\em linearly mapped coercivity}, 
providing a foundation for proving the existence of weak solutions. 
Our approach departs from traditional methods by avoiding the standard decomposition 
into a divergence-free velocity equation and a subsequent pressure equation.

Next, we address nonlinear saddle-point problems, incorporating a nonlinear mapping $\Phi$
to extend coercivity to {\em mapped coercivity}. We apply this generalized concept to prove the existence 
of solutions for various finite element discretizations of the Navier-Stokes system, including pressure-stabilization 
techniques by Brezzi-Pitk\"aranta \cite{BrezziPitkaranta1984}, 
Becker-Hansbo \cite{BeckerHansbo2008}, and local projection stabilization (LPS) 
\cite{BeckerBraack2001}. For cases involving 
dominant convection, additional velocity stabilization via LPS terms is considered \cite{BraackErn2006}. 
We establish the existence of solutions for arbitrary data in a concise manner.

The structure of this paper is as follows: Section \ref{sec:mappedcoercivity}
provides an overview of key results from \cite{BeckerBraack2024}, 
forming the foundation for our subsequent analysis. Section \ref{sec:linearsaddlepoint}.
applies this methodology to linear saddle-point problems, exploring the connection between the inf-sup condition 
and {\em linearly mapped coercivity}. As a case study, we analyze the Stokes equation and 
several equal-order finite element schemes. Section \ref{sec:nonlinearsaddlepoints}
then extends this framework to nonlinear saddle-point problems. 
For nonlinearities of at most quadratic order, we construct a surjective mapping $\Phi: X \to X$
to ensure {\em mapped coercivity}, applying this approach to the Navier-Stokes system. 
We further extend the analysis to cover a range of stabilized finite element methods, including the aforementioned pressure-stabilization methods combined with convection stabilization via LPS.

\section{Mapped coercivity for nonlinear operators}\label{sec:mappedcoercivity} 

Let $X$ and $Y$ be two separable Banach spaces with their respective dual spaces $X'$ and $Y'$. The space $X$ is assumed to be reflexive, with a strictly convex dual
and with a Schauder basis.
Let $A :X\to Y'$ be the operator under investigation. For given $b\in Y'$, we are interested in solutions
of
\beqal\label{eq:1}
   x\in X: \quad A(x) &=& b.
\eeqal
In the following we use two types of mapped coercivity, see \cite{BeckerBraack2024}:
\begin{definition}
\label{def:linearlymappedcoercivity}
An operator $A:X\to Y$ is called {\em linearly mapped coercive}, if there exists
$\Phi\in L(X,Y)$ which is onto and satisfies
\beqal\label{eq:linearlymappedcoercivity}
  \liminf_{x\in X, \normX{x}\to \infty}\frac{\langle A(x),\Phi x\rangle_{Y',Y}}{\normX{x}} &=& +\infty.
\eeqal
\end{definition}
The second type of coercivity allows nonlinear maps $\Phi:X\to Y$ which are of a particular form.
Due to the strict convexity of $X'$, there are projections $P_n:X\to X_n$ and a duality map $\Psi:X\to X'$ with the following properties
\beqal\label{eq:psibounded}
  \normXS{\Psi(x)} &\le & C\normX{x},\\
\label{eq:psicoercive}
   \langle \Psi(x),x\rangle_{X',X} &>& 0\quad\forall x\in X\setminus\{0\},\\
   \langle \Psi(z),x-P_nx\rangle_{X',X} &=& 0\quad\forall z\in X_n,\ \forall x\in X.\label{eq:psiprojected}
\eeqal
We define the following property:
\begin{definition}
\label{def:mappedcoercivity}
An operator $A:X\to Y$ is called {\em mapped coercive}, if there exists an operator
$\Psi:X\to X'$ satisfying (\ref{eq:psibounded})-(\ref{eq:psiprojected})  and a surjection
$S\in\mathcal L(X',Y)$ s.t.:
\beqal\label{eq:mappedcoercivity}
  \liminf_{x\in X, \normX{x}\to \infty}\frac{\langle A(x),S\Psi(x)\rangle_{Y',Y}}{\normX{x}} &=& +\infty.
\eeqal
\end{definition}
The property of mapped coercivity in combination with weak-*-continuity of $A$
\beqal\label{eq:weakcont}
   x_n\rightharpoonup x \mbox{ in }X &\Longrightarrow&  A(x_{n}) \stackrel{*}{\rightharpoonup} A(x)\mbox{ in }Y'
\eeqal
is sufficient to get the existence of solutions of (\ref{eq:1}), see \cite{BeckerBraack2024}:
\begin{theorem}[Existence of solutions]\label{thm:1}
Let $X$ and $Y$ be two separable real Ba\-nach spaces, $X$ reflexive. 
$A:X\to Y'$ is assumed to be weak-$*$-continuous (\ref{eq:weakcont}) and linearly mapped coercive
(Def. \ref{def:linearlymappedcoercivity}) or mapped coercive 
(Def. \ref{def:mappedcoercivity}). Then, for any $b\in Y'$, equation (\ref{eq:1}) has  a solution
$x\in X$.
\end{theorem}
In the next lemma, we treat the particular case of real, separable Hilbert space $X$ with
inner product $(\cdot,\cdot)_X$ and $Y=X$. An important role
has the Riesz isomorphism $\Psi_{iso}\in L(X,X')$, defined by
\beqa
   \langle\Psi_{iso} u,v\rangle_{X',X} &:=& (u,v)_X\qquad\forall v\in X.
\eeqa
\begin{lemma}\label{lemma:equivcoercive}
Let $X$ be a separable Hilbert spaces and $A:X\to X'$. Then the following two properties are equivalent:
\begin{description}
\item[(a)] $A$ is linearly mapped coercive.
\item[(b)] $A$ is mapped coercive with $\Psi=\Psi_{iso}\in L(X,X')$ the Riesz isomorphism.
\end{description}
\end{lemma}
\begin{proof}
First of all, we like to emphasize that the injective Riesz isomorphism $\Psi_{iso}\in L(X,X')$ satisfies
 (\ref{eq:psibounded})-(\ref{eq:psiprojected}). \\
 $(a)\Rightarrow (b)$: The property of mapped coercivity follows with the setting $\Psi=\Psi_{iso}$ and $S:=\Phi\circ\Psi_{iso}^{-1}$,
 because it implies $S\circ\Psi =\Phi$.\\
$(b)\Rightarrow (a)$: The property of linearly mapped coercivity follows with the setting $\Phi:=S\circ\Psi_{iso}$: Since $\Psi_{iso}$ 
is an isometry and $S$ in onto, $\Phi$ is onto as well.
\end{proof}

However, in the case of Hilbert spaces, this lemma does not imply that the two properties in
Definition  \ref{def:linearlymappedcoercivity} and  \ref{def:mappedcoercivity} are equivalent.
In particular, the mapping $\Psi$ must not necessarily be the Riesz isomorphism.
Also a nonlinear scaling of $\Psi_{iso}$ is possible.
For any strictly positive, bounded function $\sigma:X\to\nR_{>0}$,
 $0<\sigma(x)\le M$ for all $x\in X$, the product function $\Psi:X\to X'$, defined by
 $\Psi(x):=\sigma(x)\Psi_{iso}(x)$, also   satisfies
 (\ref{eq:psibounded})-(\ref{eq:psiprojected}). 

The following theorem states the existence of solutions leveraging the property of being mapped coercive.
The result can be interpreted as having existence for bounded (or small) data.

For a projection $P_n\in L(X,X_n)$ and $S\in L(X',Y)$ let $R_n:X_n\to X_n$ given by 
\beqal\label{eq:Rn}
   R_n(x) &:=& P_nS^*(A(x)-b).
\eeqal
Let us consider the discrete nonlinear eigenvalue problem
\beqal\label{eq:ev}
   x\in X_n: \qquad R_n(x)&=&\lambda x.
\eeqal
for $\lambda\in\nR$.
The following lemma makes a statement about the boundedness of such eigen-solutions.
\begin{lemma}\label{lemma:ev}
Let $A:X\to Y$ and $b\in Y'$.  We assume the existence of an operator
$\Psi:X\to X'$ satisfying (\ref{eq:psibounded}) and (\ref{eq:psiprojected}), 
$S\in\mathcal L(X',Y)$ and $r>0$ s.t.
for any $x\in X$ with $\normX{x}\ge r$ it holds
\beqa
   \langle A(x)-b,S\Psi(x) \rangle_{Y',Y} &\ge& 0.
\eeqa
Then, any solution of (\ref{eq:ev}) with $\lambda<0$ is bounded $\normX{x}<r$.
\end{lemma}
\begin{proof}
The proof is done by contradiction. Let us assume the existence of a solution of  (\ref{eq:ev}) with $\lambda<0$ and $\normX{x}\ge r$.
We obtain the contradiction
\beqa
   0&>&\lambda \langle \Psi(x),x\rangle_{X'\times X}\\
   &=&  \langle \Psi(x),R_n(x)\rangle_{X'\times X}\\
   &=&  \langle A(x)-b,S\Psi(x)\rangle_{Y',Y}
   \ \ge\ 0.
\eeqa
\end{proof}
\def\dpXn#1#2{\langle#1,#2\rangle_{X_n',X_n}}
\def\dpX#1#2{\langle#1,#2\rangle_{X',X}}
\def\dpY#1#2{\langle#1,#2\rangle_{Y',Y}}

Now, we consider the discrete equations
\beqal\label{eq:discrete}
   x\in X_n:\quad \dpY{A(x)}{y} = \dpY{b}{y}\quad\forall y\in Y_n,
\eeqal
with
\beqa
   Y_n &:=& (S\circ P_n^*)(X_n').
\eeqa
\begin{lemma}\label{lemma:11}
Let $S\in\mathcal L(X',Y)$, and $R_n: X_n\to X_n$ given by (\ref{eq:Rn}).
Then, $x\in X_n$ is solution of  (\ref{eq:discrete})  iff $R_n(x)=0$.
\end{lemma}
\begin{proof}
We have the following equivalences:
\beqa
   R_n(x)=0 &\Longleftrightarrow&  P_nS^*(A(x)-b)=0\\
    &\Longleftrightarrow&  \dpXn{x'}{ P_nS^*(A(x)-b)} = 0\quad\forall x'\in X_n'\\
    &\Longleftrightarrow&  \dpY{A(x)-b}{SP_n^*x'} = 0\quad\forall x'\in X_n'\\
    &\Longleftrightarrow&  \dpY{A(x)-b}{y} = 0\quad\forall y\in Y_n.
\eeqa
This completes the proof.
\end{proof}
\begin{theorem}\label{thm:almostcoercivity}
Let $A:X\to Y$ be weak-*-continuous and $b\in Y'$.  We assume the existence of an operator
$\Psi:X\to X'$ satisfying (\ref{eq:psibounded})-(\ref{eq:psiprojected}), a surjection
$S\in\mathcal L(X',Y)$ and $r>0$ s.t.
for any $x\in X$ with $\normX{x}\ge r$ it holds
\beqa
   \langle A(x)-b,S\Psi(x) \rangle_{Y',Y} &\ge& 0.
\eeqa
Then there exists a solution $x\in X$ of $A(x)=b$ with $\normX{x}\le r$.
\end{theorem}
\begin{proof}
(a)  We firstly show the existence of bounded solutions $x_n\in X_n$, $\normX{x_n}\le r$, of the 
finite dimensional problems (\ref{eq:discrete}).
Let us assume that no solution $x\in X_n$ of (\ref{eq:discrete}) with $\normX{x}\le r$ exists.
Then by Lemma \ref{lemma:11}, for any $x\in X_n$ with $\normX{x}\le r$ we have $R_n(x)\not=0$.
Hence, the continuous opeartor $T:X_n\to X_n$ defined by $T_n(x):=-r\normX{R_n(x)}^{-1}R_n(x)$ is well defined
for $\normX{x}\le r$. By Schauder's Fixed Point Theorem we get a fixed point $T_n(x)=x$, $\normX{x}= r$.
This fixed point solves $R_n(x)=\lambda x$ with $\lambda:=-\normX{R_n(x)}r^{-1}<0$.
By Lemma \ref{lemma:ev} we get the contradiction $\normX{x}<r$. Hence, there is a solution $x\in X_n$ of (\ref{eq:discrete}) with $\normX{x}\le r$.\\
(b) We have 
\beqa 
   \overline{\bigcup_{n\in\mathbb N} P_n^*(X_n')} &=& X.
\eeqa
Since $S$ is onto and continuous, we obtain $ \overline{\bigcup_{n\in\mathbb N} Y_n} = Y$.\\
(c) Solutions of the infinite dimensional problem follows by the boundedness of the
discrete solutions from (a), the density of the $Y_n$ from (b), and the weak-*-continuity of $A$, see \cite{BeckerBraack2024}.
\end{proof}

\section{Mapped coercivity for linear saddle-point problems}\label{sec:linearsaddlepoint} 

In this section we consider operators $A:X\to X'$ of saddle-point form. Here, $X=V\times Q$ is a product space of 
two separable Hilbert Banach spaces $V$ and $Q$. The two components
of $x:=(u,p)\in X$ are denoted by $u\in V$ and $p\in Q$. 

The operator $A$ is assumed to be of the form
\beqal\label{eq:spoa}
   Ax &=&   \begin{bmatrix}L & -B^*\\ B & 0\end{bmatrix} \begin{bmatrix}u \\ p\end{bmatrix},
\eeqal
with linear continuous operators $L\in L(V,V')$ and $B\in L(V,Q')$. 
We make the following assumptions:
\begin{description}
\item[(A1)]  $L:V\to V'$ is coercive, i.e. it exists $\alpha>0$ s.t.
\beqa
 \langle Lu,u\rangle_{V',V} &\ge& \alpha\normV{u}^2\qquad\forall u\in V.
\eeqa
\item[(A2)]  
It exists $\Theta\in L(Q, V)$ which is onto
s.t.
\beqa
   \langle B^*p,\Theta p\rangle_{V',V}\ge\normQ{p}^2.
\eeqa
\end{description}
\begin{proposition}\label{prop:saddlecoercive}
The operator $A:X\to X$ is assumed to be of saddle-point form (\ref{eq:spoa}) with the properties
(A1) and (A2). Then $A$ is linearly mapped coercive (Def. \ref{def:linearlymappedcoercivity}).
\end{proposition}
\begin{proof}
Let $\Theta \in L(Q,V)$ be the operator in (A2). 
We define
$\Phi\in L(X,X)$ for $x=(u,p)$ by $\Phi x:=(u-\tau \Theta p,p)$ with a constant $\tau>0$.
$\Phi$ is onto, because for arbitrary $x=(u,p)\in X$ it holds for $w:=u+\tau \Theta p\in V$:
\beqa
    \Phi(w,p) = \Phi(u+\tau \Theta p,p)= (u,p)=x.
\eeqa
Hence, $\Phi$ is onto.
Furthermore, 
\beqa
  \langle Ax,\Phi x\rangle_{X',X} &=&   \langle Ax,x\rangle_{X',X} - \tau\langle Lu-B^*p,\Theta p\rangle_{V',V}\\
   &=& \langle Lu,u\rangle_{V',V}
   -\tau\langle Lu,\Theta p\rangle_{V',V} +\tau\langle B^*p,\Theta p\rangle_{V',V} \\
   &\ge&  \alpha\normV{u}^2
   +\tau\normQ{p}^2-\tau c_{\Theta }\normVS{Lu}\normQ{p}.
 \eeqa
Young's inequality yields $\tau c_{\Theta }\normVS{Lu}\normQ{p}\le \frac\alpha2\normV{u}^2+\frac1{2\alpha}\tau^2 c_L^2c_{\Theta }^2\normQ{p}^2$ and, hence,
\beqa
  \langle Ax,\Phi x\rangle_{X',X} &\ge&  
  \frac\alpha2\normV{u}^2
  +\tau\left(1-\frac{\tau}{2\alpha}c_{\Theta }^2c_L^2\right)\normQ{p}^2.
\eeqa
For $\tau>0$ sufficiently small we get with $\beta:=\min(\alpha,\tau)/2$
\beqa
   \langle Ax,\Phi x\rangle_{X',X}   &\ge&  \frac\alpha2\normV{u}^2+\frac{\tau}2\normQ{p}^2\ \ge\ \beta\normX{x}^2.
\eeqa
This shows that $A$ is linearly mapped coercive (\ref{eq:linearlymappedcoercivity}).
\end{proof}
\begin{corollary}\label{cor:linearsaddle}
The operator $A:X\to X$ is assumed to be of saddle-point form (\ref{eq:spoa}) with the properties
(A1) and (A2). Then $A$ is linearly mapped coercive (Def. \ref{def:linearlymappedcoercivity}).
\end{corollary}
\begin{proof}
The assertion follows from Proposition \ref{prop:saddlecoercive} and Lemma \ref{lemma:equivcoercive}.
\end{proof}

\subsection{\bf Relation to linear saddle-point problems with inf-sup condition}\label{sec:infsup} 

The following Lemma relates the condition (A2) from above to the well-known inf-sup condition.
\begin{lemma}\label{lemma:infsup}
The following three properties are equivalent for any $B\in L(V,Q')$:
\begin{enumerate}
\item[(a)]  (A2) is valid.
\item[(b)] $B$ satisfies an inf-sup condition, i.e. it exists $\gamma>0$ s.t.
\beqal\label{eq:infsup}
 \inf_{p\in Q\setminus\{0\}}\sup_{v\in V\setminus\{0\}} \frac{\langle Bv,p\rangle_{Q',Q}}{\normV{v}\normQ{p}} &\ge& \gamma.
\eeqal
\item[(c)] $B$ is an isomorphism from $(Ker\,B)^\perp$ to $Q'$ and $\normV{v}\le c_B\normQS{Bv}$ for all $v\in (Ker\, B)^\perp$ with a constant $c_B\ge0$.
\end{enumerate}
\end{lemma}
\begin{proof}
(a)$\Longrightarrow$(b): Let $p\in Q$ be given. We define $v:=\Theta p\in V$. By continuity of $\Theta $ we have
$\normV{v}\le c_\Phi\normQ{p}$. With (A2) we derive
\beqal\label{eq:Bvp}
   \langle Bv,p\rangle_{Q',Q} &=&  \langle B^*p,\Theta p\rangle_{Q',Q}  \ge \normQ{p}^2\ge \frac1{c_\Phi}\normQ{p}\normV{v}.
\eeqal
This implies the  inf-sup condition with $\gamma:=c_\Phi^{-1}$. \\
 (b)$\Longleftrightarrow$(c): For this equivalence, we refer to \cite{GiraultRaviart1986}.\\
 (c)$\Longrightarrow$(a): Let us denote the Riesz Isomorphism on $Q$ by $R:Q\to Q'$, and the restriction of $B$
 on $(Ker\,B)^\perp$ by $B_0$. With (c) we now that $B_0\in L((Ker\,B)^\perp,Q')$ is an isomorphism. We
 now define $\Theta:=B_0^{-1}\circ R\in L(Q,V)$. We obtain
 \beqa
      \langle B^*p,\Theta p\rangle_{V',V} &=&   \langle (B\circ\Theta )p,p\rangle_{Q',Q} \\
      &=&   \langle (B\circ B_0^{-1}\circ R)p,p\rangle_{Q',Q}\\
      &=&   \langle Rp,p\rangle_{Q',Q}\ =\ \normQ{p}^2.
 \eeqa
\end{proof}

\subsection{\bf Generalized Stokes problem}\label{sec:stokes} 

As an application, we may consider the Stokes problem
\beqal\label{eq:stokes1}
  -\Delta u+\nabla p &=& f\quad\mbox{in }\Omega,\\
   \mbox{div}\,u &=& g\quad\mbox{in }\Omega,
   \\
   u&=& 0\quad\mbox{on }\partial\Omega.
\eeqal
The corresponding space is $X:=H^1_0(\Omega)\times L^2_0(\Omega)$, and
$A:X\to X$ is of saddle-point form (\ref{eq:spoa}) with 
$L:=-\Delta$ and $B:={\rm div}$.
$L$ is coercive, hence (A1) holds.
The inf-sup condition (\ref{eq:infsup}) is fulfilled. 
Hence, according to Lemma \ref{lemma:infsup}, (A2) holds. By Corollary \ref{cor:linearsaddle},
for any $f\in H^{-1}(\Omega)$, $g\in L^2(\Omega)$, we get existence of a solution $(u,p)\in X$.

\subsection{\bf Stabilized linear saddle-point problems}\label{sec:generalizedspp} 

In some cases, condition (A2) does not hold for subspaces $X_h\subset X$, although valid in $X$.
Therefore, numerical schemes for solving a saddle point problem in a (usually finite dimensional) subspace are often based on 
introducing a stabilized operator. Let $A\in L(X,X')$  of the form (\ref{eq:spoa}) and $A_h\in L(X_h,X_h')$ of the form
\beqal\label{eq:spogeneralized0}
   \langle A_hx,y\rangle_{X_h,X_h'} :=   \langle Ax,y\rangle_{X,X'}+ \langle Tp,q\rangle_{X_h,X_h'}
\eeqal
for $x=(u,p), y=(v,q)\in X_h:=V_h\times Q_h$. Hence, $A_h$ can be expressed as 
\beqal\label{eq:spogeneralized}
   A_hx &=&   \begin{bmatrix}L & -B^*\\ B & T\end{bmatrix} \begin{bmatrix}u \\ p\end{bmatrix},
\eeqal
with linear continuous operators $L\in L(V,V')$, $B\in L(V,Q')$, as before,  and $T\in L(Q_h,Q_h')$. 
Making the following additional assumption:
\begin{description}
\item[(A2a)]  $ \langle Tp,p\rangle_{Q_h',Q_h}\ge 0$,
\end{description}
then, we will relax condition (A2)  to
\begin{description}
\item[(A2b)]  
It exists $\Theta\in L(Q_h, V_h)$ which is onto 
and a constant $c_T\in\nR$ s.t.
\beqa
   \langle B^*p,\Theta p\rangle_{V',V}\ge\frac12\normQ{p}^2-c_T\langle Tp,p\rangle_{Q_h',Q_h}\qquad\forall p\in Q_h.
\eeqa
\end{description}
Obviously, in the case $X_h=X$, (A2) is equivalent to (A2a) in combination with (A2b) and $T=0$.\\
\begin{corollary}\label{cor:spogeneralized}
The operator  $A_h\in L(X_h,X_h')$ is assumed to be of saddle-point form (\ref{eq:spogeneralized}) with the properties
(A1), (A2a) and (A2b).
Then $A_h$ is linearly mapped coercive (Def. \ref{def:linearlymappedcoercivity}).
\end{corollary}
\begin{proof}
The choice of $\Phi$ is the same as the proof of Proposition \ref{prop:saddlecoercive} (i.e. the linear case),
 $\Phi x:=(u-\tau\Theta p,p)$. With the technique in the proof of Proposition \ref{prop:saddlecoercive} 
 and with use of (A2a) and (A2b): 
 \beqa
 \langle A_hx,\Phi x\rangle_{X',X} &=&   \langle Lu,u\rangle_{V',V}
   -\tau\langle Lu,\Theta p\rangle_{V',V}\\
   && +\tau\langle B^*p,\Theta p\rangle_{V',V}+\langle Tp,p\rangle_{Q_h',Q_h} \\
   &\ge&  \frac\alpha2\normV{u}^2+\tau\left(1-\frac{\tau}{2\alpha}c_{\Theta }^2c_L^2\right)\normQ{p}^2 +  (1-\tau c_T)\langle Tp,p\rangle_{Q_h',Q_h}.
\eeqa
For $0<\tau\le \min(c_T^{-1},\alpha(c_\Theta c_L)^{-2})$ we end up with linearly mapped coercivity (\ref{eq:linearlymappedcoercivity}):
\beqa
  \langle A_hx,\Phi x\rangle_{X',X} 
   &\ge&   \frac\alpha4\normV{u}^2+ \frac\tau2\normQ{p}^2\qquad\forall x\in X_h.
\eeqa
\end{proof}
\begin{lemma}\label{lemma:a2binterpol}
We assume (A2) for $B\in L(V,Q')$. Further, let $T\in L(Q_h,Q_h')$ be given,
the existence of a projection $\Pi_h\in L(V,V_h)$ and a constant $c_I\ge0$ s.t.
\beqal\label{eq:a2binterpol}\qquad
   \langle B^*p,\Pi_hu-u\rangle_{V',V} &\le& c_I \langle Tp,p\rangle_{Q_h',Q_h}^{1/2} \normV{u}\qquad\forall (u,p)\in V\times Q_h.
\eeqal
Then $B\in L(V_h,Q_h')$ fulfills (A2b).
\end{lemma}
\begin{proof}
Due to (A2) for the space $X$, we have the existence of $\Theta \in L(Q,V)$ such that $\langle B^*p,\Theta  p\rangle_{V',V}\ge \normQ{p}^2$ for all $p\in Q_h$.
We define  $\Theta_h:=\Pi_h\circ \Theta|_{Q_h} \in L(Q_h,V_h)$.  Since $\Theta$ is onto and $\Pi_h$ is projection, $\Theta_h$ is onto.
Further, we get for any $p\in Q_h$
by triangle inequality and (\ref{eq:a2binterpol}):
\beqa
   \langle B^*p,\Theta_h p\rangle_{V',V} &= &  
    \langle B^*p,\Theta p\rangle_{V',V} + \langle B^*p,(\Pi_h-id)\Theta p\rangle_{V',V} \\
   &\ge&   \normQ{p}^2- c_I \langle Tp,p\rangle_{Q_h',Q_h}^{1/2} \normV{\Theta p}.
\eeqa
With the continuity of $\Theta$ and Young's inequality we arrive at
\beqa
    \langle B^*p,\Theta_h p\rangle_{V',V} &\ge &  \frac12\normQ{p}^2-c_T \langle Tp,p\rangle_{Q_h',Q_h}.
\eeqa 
\end{proof}
\begin{corollary}\label{cor:ahlinearmapped}
$A\in L(X,X')$ is assumed to be of saddle-point form (\ref{eq:spoa}) satisfying properties (A1) and (A2).
$A_h\in L(X_h,X_h')$ is assumed to be of the form (\ref{eq:spogeneralized0}) with the properties
(A2a) and (\ref{eq:a2binterpol}).
Then $A_h$ is linearly mapped coercive (Def. \ref{def:linearlymappedcoercivity}).
\end{corollary}
\begin{proof}
The assertion follows directly with Lemma \ref{lemma:a2binterpol} and Corollary \ref{cor:spogeneralized}.
\end{proof}

\subsection{\bf Equal-order finite element discretization of the generalized Stokes problem}\label{sec:stokesfem} 

We now consider conforming finite elements, i.e. finite dimensional
subspaces $X_h:=V_h\times Q_h\subset X=H^1_0(\Omega)\times L^2_0(\Omega)$ consisting of piece-wise polynomials
of degree $k$. 
The parameter $h$ stands for the maximal mesh size of a triangulation $\mathcal T_h$
of the polygonal bounded domain $\Omega$
in cells $T\in\mathcal T_h$ and outer radius $h_T:=\mbox{diam}(T)$. 
All cells $T$ may consist of triangles, quadrilaterals, tetrahedra or hexahedra.
Because not all subspaces, in particular equal-order finite elements,
allow a discrete analogon of the inf-sup property, we need a discrete version of the mapping $\Phi$.

The operator $A_h$ is of the form (\ref{eq:spogeneralized}) with operators $L\in L(V,V')$, $B\in L(V,Q')$, and  $T\in L(Q_h,Q_h')$ 
as a stabilization term, because the divergence operator
$B$ is not discrete inf-sup stable for equal order finite elements.
In the following subsections, we consider different types of stabilization.
With Corollary \ref{cor:ahlinearmapped} we then get the existence of a discrete solution for the corresponding
finite element approximation of the generalized Stokes problem.
\begin{proposition}\label{prop:pspg}
The generalized Stokes equation (\ref{eq:stokes1}) discretized by equal-order $P_1$-ele\-ments and stabilized by
method with the property
\beqal\label{eq:Tpp}
   \langle Tp,p\rangle_{Q_h',Q_h} &\ge c\sum_{K\in \mathcal T_h} h_K^2\|\nabla p\|_{L^2(K)}^2\qquad\forall p\in Q_h,
\eeqal
with an arbitrary positive constant $c>0$,
results into a linearly mapped coercive operator $A_h$.
\end{proposition}
\begin{proof}
The operator $T\in L(Q_h,Q_h')$ obviously satisfies (A2a). 
In order to validate condition (A2b) we will use Corollary \ref{cor:ahlinearmapped}.
Let $\Pi_h:V\to V_h$ the $L^2-$projection. Due to the cell-wise $L^2$-estimate for $K\in\mathcal T_h$,
\beqa
   \| (\Pi_h-id)u\|_{L^2(K)} &\le& c_I h_K\|\nabla u\|_{L^2(K)}\quad\forall u\in H^1(K),
\eeqa
we have the upper bound
\beqa
   \langle B^*p,(\Pi_h-id)u)\rangle_{V',V} &= &\sum_{K\in \mathcal T_h} \int_K \nabla p  (\Pi_h-id)u\\
    &\le &\sum_{K\in \mathcal T_h} \| \nabla p\|_{L^2(K)} \|(\Pi_h-id)u\|_{L^2(K)}\\
   &\le & c_I\sum_{K\in \mathcal T_h} \|\nabla p\|_{L^2(K)}h_K \| \nabla u\|_{L^2(K)}\\
   &\le & c_I\left(\sum_{K\in \mathcal T_h} h_K^2\|\nabla p\|_{L^2(K)}^2\right)^{1/2} \normV{u}.
\eeqa
With (\ref{eq:Tpp}) we arrive at (\ref{eq:a2binterpol}), and Corollary \ref{cor:ahlinearmapped} ensures the linearly mapped coercivity.
\end{proof}

\subsubsection{\bf Brezzi-Pitkar\"anta pressure stabilization}\label{sec:bps}

A popular pressure stabilization \cite{BrezziPitkaranta1984} of the Stokes equations with piece-wise linear ($k=1$) 
equal-order finite elements uses for $p,q\in Q_h$ the stabilization term
\beqal\label{eq:pspg}
   \langle Tp,q\rangle_{Q_h',Q_h} &:=& \sum_{K\in \mathcal T_h} h_K^2\int_K \nabla p\cdot\nabla q.
\eeqal
\begin{corollary}\label{cor:pspg}
The generalized Stokes equation (\ref{eq:stokes1}) discretized by equal-order $P_1$-ele\-ments and stabilized by
method (\ref{eq:pspg}) satisfies conditions (A1), (A2a) and (A2b). Hence, the corresponding operator $A$ is linearly
mapped coercive.
\end{corollary}
\begin{proof}
The assertion follows directly from Proposition \ref{prop:pspg} and the particular form of $T$ in (\ref{eq:pspg}).
\end{proof}

\subsubsection{\bf A simple stabilization method for Stokes}\label{sec:bhs}

The pressure-velocity stabilization in \cite{BeckerHansbo2008} for $P_1$ elements reads
\beqal\label{eq:bhs}
    \langle Tp,q\rangle_{Q_h',Q_h} &:= & \int_\Omega \left(I_h^1(pq)-pq\right).
\eeqal
with the Lagrange interpolation operator $I_h^1:Q\to Q_h$.
\begin{corollary}
The linear operator $A_h$ corresponding to the generalized Stokes problem (\ref{eq:stokes1})
discretized with  $P_1$ equal-order elements and stabilization (\ref{eq:bhs})  
is linearly mapped coercive.
\end{corollary}
\begin{proof}
Proposition 1 in \cite{BeckerHansbo2008} gives for  any cell $K\in\mathcal T_h$ and any linear function $p$ the estimate
\beqa
     h_K^2\|\nabla p\|_{L^2(K)}^2  &\le& c\int_K \left(I_h^1(p^2)-p^2\right).
\eeqa
The assertion now follows directly from Proposition \ref{prop:pspg}.
\end{proof}

\subsubsection{\bf Local projection stabilization (LPS)}\label{sec:lps}

Another method which can be cast into this framework is the local pressure stabilization (LPS), see \cite{BeckerBraack2001}. 
Let $k$ be the polynomial order of the discrete pressure and velocity, i.e.
equal-order elements of order $k$.
We consider the space $D_h$ of patch-wise polynomials of order $k-1$, discontinuous across patches of elements.
Such a patch can be considered as an element of a macro-triangulation $\mathcal T_{2h}$. 
Let $\pi_h:L^2(\Omega)^d\to D_h$ the $L^2$-projection.
The LPS method makes use of the stabilization term
\beqal\label{eq:lps}
    \langle Tp,q\rangle_{Q_h',Q_h} &:= & \sum_{K\in \mathcal T_{2h}}h_K^2\int_K  \left[\nabla p-\pi_h(\nabla p)\right]\cdot \left[\nabla q-\pi_h(\nabla q)\right].
\eeqal
\begin{proposition}
The linear operator corresponding to the generalized Stokes problem (\ref{eq:stokes1})
discretized with  equal-order elements of order $k$ and LPS (\ref{eq:lps})  
satisfies conditions (A1), (A2a) and (A2b). Hence, the corresponding operator $A$ is linearly
mapped coercive.
\end{proposition}
\begin{proof}
(A1) and (A3) are obviously satisfied. Condition (A2b) has still to be validated with $T$ given in (\ref{eq:lps}).
We use the
projection operator $j_h:V\to V_h$ from \cite{BeckerBraack2001} with the orthogonality property
\beqal\label{eq;jhortho}
   \int_\Omega d_h (j_hv-v)  &=& 0\quad\forall d_h\in D_h,\ \forall v\in V.
\eeqal
With the orthogonality property (\ref{eq;jhortho}) and a standard interpolation property
of $j_h$ (see e.g. \cite{BeckerBraack2001}) and due to $\pi_h(\nabla p_h)\in D_h$ we get 
\beqa
   \langle B^*p,(j_h-id)u)\rangle_{V',V}&= &\sum_{K\in \mathcal T_{2h}} \int_K \nabla p  (j_h-id)u\\
   &=& \sum_{K\in \mathcal T_{2h}} \int_K (\nabla p-\pi_h(\nabla p))  (j_h-id)u\\
    &\le &\sum_{K\in \mathcal T_{2h}} \| \nabla p-\pi_h(\nabla p)\|_{L^2(K)} \|(j_h-id)u\|_{L^2(K)}\\
   &\le & c_I\sum_{K\in \mathcal T_{2h}} \|\nabla p-\pi_h(\nabla p)\|_{L^2(K)}h_K \| \nabla u\|_{L^2(K)}\\
   &\le & c_I\left(\sum_{K\in \mathcal T_{2h}} h_T^2\|\nabla p-\pi_h(\nabla p)\|_{L^2(K)}^2\right)^{1/2} \normV{u}\\
   &\le & c_I  \delta_0^{-1/2}\langle Tp,p\rangle_{Q_h',Q_h}^{1/2} \normV{u}.
\eeqa
With Lemma \ref{lemma:a2binterpol} we get (A2b).
\end{proof}

\section{Mapped coercivity for nonlinear saddle-point problems}\label{sec:nonlinearsaddlepoints} 

In this subsection we consider operators $A:X\to X'$  of the form
\beqal\label{eq:spoanonlinear}
   A(x) &=&  \left(\begin{array}{c}C(u) \\
   			0\end{array}\right) +  \begin{bmatrix}L & -B^*\\ B & T\end{bmatrix} \begin{bmatrix}u \\ p\end{bmatrix},
\eeqal
with linear continuous operators $L\in L(V,V')$, $B\in L(V,Q')$, as before,  $T\in L(Q,Q')$ and a (possible) nonlinear
operator $C:V\to V'$. 
Let the right hand side of the form $b=(f,0)$ with $f\in V'$.
This corresponds to a problem of seeking $u\in V$ and $p\in Q$ s.t.
\beqa
   C(u)+Lu-B^*p &=& f,\\
   Bu + Tp&=& 0.
\eeqa
We make the following additional assumption:
\begin{description}
\item[(A3)]  $ \langle C(u),u\rangle_{V',V}\ge 0$.
\end{description}

\subsection{\bf Saddle-point problems with moderate non-linearity}\label{sec:moderate} 

The following Corollary allows for the nonlinear operator $C$ only linear growth:
\begin{corollary}\label{cor:moderate}
The operator $A:X\to X$ is assumed to be of saddle-point form (\ref{eq:spoanonlinear}) with the properties
(A1), (A2a), (A2b), (A3), and there is a constant  $c_N\ge 0$ s.t.
\beqal\label{eq:Clinearbounded}
  \|C(u)\|_{V'} &\le&c_N \normV{u}\qquad\forall u\in V.
\eeqal
Then $A$ is linearly mapped coercive (Def. \ref{def:linearlymappedcoercivity}).
\end{corollary}
\begin{proof}
The choice of $\Phi$ is the same as the proof of Proposition \ref{prop:saddlecoercive} (i.e. the linear case),
 $\Phi(x):=(u-\tau\Theta p,p)$. With the technique in the proof of Proposition \ref{prop:saddlecoercive} 
 and with use of (A3) and (A2b) for $\tau>0$ sufficiently small: 
\beqa
  \langle A(x),\Phi x\rangle_{X',X} 
   &\ge&  \frac\alpha2\normV{u}^2+\frac{\tau}2\normQ{p}^2+\langle C(u),u- \tau\Theta p\rangle_{V',V} \\
  && +  (1-\tau c_T)\langle Tp,p\rangle_{Q',Q}\\
   &\ge&   \frac\alpha2\normV{u}^2+• \frac{\tau}2\normQ{p}^2+  (1-\tau c_T)\langle Tp,p\rangle_{Q',Q}\\
    &&- \tau c_N\normV{u}  c_\Theta\normQ{ p}\\
   &\ge&   \frac\alpha4\normV{u}^2+ \tau\left(\frac12-\frac1\alpha \tau c_N^2c_\Theta^2\right)\normQ{p}^2\\
   &&+  (1-\tau c_T)\langle Tp,p\rangle_{Q',Q}
\eeqa
With a possible further diminuition of $\tau>0$ we end up with linearly mapped coercivity (\ref{eq:linearlymappedcoercivity}):
\beqa
  \langle A(x),\Phi x\rangle_{X',X} 
   &\ge&   \frac\alpha4\normV{u}^2+ \frac\tau4\normQ{p}^2.
\eeqa
\end{proof}

As an example we may choose the following nonlinear variant of an Oseen problem 
\beqa
  (b(u)\cdot\nabla) u-\Delta u+\nabla p &=& f\quad\mbox{in }\Omega,\\
   \mbox{div}\,u &=& 0\quad\mbox{in }\Omega,
   \\
   u&=& 0\quad\mbox{on }\partial\Omega,
\eeqa
where $b:H^1_0(\Omega)^d\to H^1(\Omega)^d$ is a vector field with the property $\mbox{div}\,b(u)=0$ to ensure (A3), and
$\|b(u)\|_{L^2(\Omega)}\le C$ for any $u\in H^1_0(\Omega)^d$ to ensure (\ref{eq:Clinearbounded}).
We obtain  for all $f\in H^{-1}(\Omega)^d$  the existence of solutions by Corollary \ref{cor:moderate} and Theorem \ref{thm:1}.

\subsection{\bf Saddle-point problems with quadratic non-linearity}\label{sec:quadrativ} 

If  (\ref{eq:Clinearbounded}) does not hold, i.e. $C(u)$ increases more than linearly, the proof of the previous proposition must be modified.
Instead of a linear map $\Phi$ we will construct a nonlinear map $S\circ\Psi:X\to Y$ as a composition of two operators to obtain
mapped coercivity according to Definition \ref{def:mappedcoercivity}. In particular,
we will choose $\tau$ in the proof of the previous Corollary as a continuous functional depending in a nonlinear way of $\normV{u}$
s.t. $\lim_{\normV{u}\to\infty}\tau=0$. Since $S$ must remain linear and $\Psi$ is allowed to be nonlinear in Definition 
\ref{def:mappedcoercivity}, we must merge this parameter into the map $\Psi$.

\begin{lemma}\label{lemma:SPsi}
Let  $X=V\times Q$, $\Theta \in L(Q,V)$ and
$\tau:V\to (0,(2c_\Theta)^{-1}]$ be a given continuous function. Then it exists a surjection $S\in L(X',X)$
and $\Psi:X\to X'$ with the properties (\ref{eq:psibounded})-(\ref{eq:psiprojected}) s.t.
\beqa
  (S\circ\Psi) (x) &=& (u - \tau(u)\Theta p,p).
\eeqa
\end{lemma}
\begin{proof}
Let us now denote the Riesz isomorphisms in $V$ and $Q$ by
$\psi_V\in L(V,V')$ and $\psi_Q\in L(Q,W')$, respectively. 
We define 
\beqal\label{eq:defpsinse}
    \Psi((u,p)) &:= & (\psi_V(u-\tau(u)\Theta p),\psi_Qp)
 \eeqal
satisfying (\ref{eq:psibounded}) and (\ref{eq:psiprojected}). In order to verify  (\ref{eq:psicoercive})
we estimate
\beqa
  \langle\Psi(x),x\rangle_{X', X} &=& \langle\psi_Vu-\tau(u) (\psi_V\circ\Theta) p,u\rangle_{V', V} + \langle\psi_Qp,p\rangle_{Q', Q} \\
  &=& \normV{u}^2-\tau(u)\langle(\psi_V\circ\Theta)p ,u\rangle_{V', V} + \normQ{p}^2\\
  &\ge&  \normV{u}^2-\tau(u) c_{\Theta }\normQ{p}\normV{u}+ \normQ{p}^2\\
  &\ge& \tfrac12\left(1-c_{\Theta } \tau(u)\right)\left(\normV{u}^2+ \normQ{p}^2\right).
\eeqa
With the bound $\tau(u)\le (2c_{\Theta })^{-1}$ we arrive at $\langle\Psi(x),x\rangle_{X', X}\ge\frac14 \normX{x}^2$
and  (\ref{eq:psicoercive}).
Moreover, $S:=\psi_V^{-1}\times \psi_Q^{-1}$ is an 
isomorphism, and in particular, it is is onto. We have
\beqa
   (S\circ\Psi) (x) &=& (\psi_V^{-1}(\psi_V(u-\tau(u)\Theta p)),\psi_Q^{-1} \psi_Qp) \
   =\ (u - \tau(u)\Theta p,p).
\eeqa
\end{proof}
\begin{proposition}\label{prop:noyetcoercive}
The operator $A:X\to X$ is assumed to be of saddle-point form (\ref{eq:spoanonlinear}) with the properties
(A1), (A2a), (A2b), (A3) and
\beqal\label{eq:Cquadraticounded}
{\bf (A4)}\qquad  \|C(u)\|_{V'} &\le&c_N \normV{u}^2\qquad\forall u\in V.
\eeqal
Let $b\in V'\times \{0\}$. Then, equation (\ref{eq:1}) has a solution. 
\end{proposition}
\begin{proof}
We will use Theorem \ref{thm:almostcoercivity} to validate the assertion.\\
{\em Step 1:} We use the composition $S\circ\Psi$ from the previous Lemma \ref{lemma:SPsi} with an appropriate scaling function 
$\tau:V\to\nR_{>0}$ to be determined below.
With Young's inequality we obtain
\beqa
   &&\langle A(x),S\Psi (x)\rangle_{X',X} \\
   &=& \langle C(u)+Lu-B^*p,u - \tau(u)\Theta p\rangle_{V',V} + \langle Bu+ Tp,p\rangle_{Q',Q}\\
   &=&\langle C(u)+Lu,u \rangle_{V',V}- \tau(u)\langle C(u)+Lu-B^*p,\Theta p\rangle_{V',V}+\langle Tp,p\rangle_{Q',Q}\\
   &\ge&  \frac\alpha2\normV{u}^2+\frac{\tau(u)}2\left(1-\frac{c_{\Theta }^2}{\alpha}\tau(u)(c_N \normV{u}+c_L)^2\right)\normQ{p}^2\\
&&   +(1-\tau(u)c_T)\langle Tp,p\rangle_{Q',Q}.
\eeqa
This gives us for any continuous $\tau:V\to (0,\tau_{max,1}]$ with 
$\tau_{max,1}:= c_T^{-1}$ and 
\beqal\label{eq:boundtau}
 \tau(u) &\le& \frac12 \alpha c_{\Theta }^{-2}(c_N \normV{u}+c_L)^{-2}
\eeqal
the lower bound
\beqal\label{eq:almostcoercivity}
  \langle A(x),S\Psi (x)\rangle_{X',X}  &\ge&  \frac\alpha2\normV{u}^2+\frac{\tau(u)}4\normQ{p}^2.
\eeqal
{\em Step 2:}
The result in Step 1 does not imply that the operator $A$ is mapped coercive, because the upper bound (\ref{eq:boundtau}) implies
$\lim_{\normV{u}\to\infty}\tau(u)=0$. 
However, we will show 
that the lower bound (\ref{eq:almostcoercivity}) implies the solvability of the equation $A(x)=b$ for our particular $b$. 
For $b=(f,0)$ with $f\in V'$ we have the upper bound
\beqa
  |\langle b,S\Psi (x)\rangle| &\le& \normVS{f}\normV{u-\tau(u)\Theta p}\\
  &\le&  \left(\frac1{\alpha}+2\tau(u)c_\Theta^2\right)\normVS{f}^2+\frac\alpha4 \normV{u}^2+\frac{\tau(u)}8\normQ{p}^2.
\eeqa
Thus, in combination with  (\ref{eq:almostcoercivity}), we obtain
\beqa
 \langle A(x)-b,S\Psi (x)\rangle 
 &\ge&  \frac\alpha4\normV{u}^2+\frac{\tau(u)}8\normQ{p}^2-\left(\frac1{\alpha}+2\tau(u)c_{\Theta }^2\right)\normVS{f}^2.
\eeqa
For any continuous $\tau:V\to (0,\tau_{max,2}]$ with $\tau_{max,2}:=\min(\tau_{max,1},\alpha^{-1} c_{\Theta }^{-2})>0$ and (\ref{eq:boundtau}) we obtain 
\beqa
 \langle A(x)-b,S\Psi (x)\rangle &\ge&  \frac\alpha4\normV{u}^2+\frac{\tau(u)}8\normQ{p}^2-\frac3{\alpha}\normVS{f}^2.
\eeqa
The choice $ \tau(u) := \min\left(\tfrac12\alpha c_{\Theta }^{-2}(c_N \normV{u}+c_L)^{-2},  \tau_{max,2}\right)$ satisfies all requirements on $\tau$.\\
{\em Step 3:} Now, we derive an appropriate choice of $r>0$ s.t. the lower bound in Theorem \ref{thm:almostcoercivity} is valid.
Let $r_1:=4\alpha^{-1}\normVS{f}$, $r_2:=r_1\alpha^{1/2}\tau_{max,3}^{-1/2}$ and $\tau_{max,3}:=\min(\alpha c_{\Theta }^{-2}(c_N r_1+c_L)^{-2}, \tau_{max,2})$.
For $\normX{x}\ge r:=\max(r_1,r_2)$ we have two cases:
\begin{description}
\item[Case 1] $\normV{u}^2\ge r_1^2/2$: We get
\beqa
   \langle A(x)-b,S\Psi (x)\rangle &\ge& \frac\alpha8r_1^2-\frac2{\alpha}\normVS{f}^2\ge  0.
\eeqa
\item[ Case 2] $\normV{u}^2<r_1^2/2$:  Then $\normQ{p}^2\ge\frac{r_2^2}2\ge\frac{\alpha r_1^2}{2\tau_{max,3}}$ and $\tau(u)\ge\tau_{max,3}$. We get
\beqa
 \langle A(x)-b,S\Psi (x)\rangle &\ge&  \frac\alpha8r_1^2-\frac2{\alpha}\normVS{f}^2\ge  0.
\eeqa
\end{description}
Hence, all requirements to apply Theorem \ref{thm:almostcoercivity}  are validated. This gives us the existence of a solution.
\end{proof}

\subsection{\bf Application to the Navier-Stokes equations}\label{sec:navierstokes} 

The stationary Navier-Stokes equation 
 in a bounded Lipschitz domain $\Omega\subset\nR^d$, $d\in\{2,3\}$,
with homogeneous Dirichlet data reads
\beqal\label{eq:nse1}
   (u\cdot\nabla) u -\mu\Delta u+\nabla p &=& f\quad\mbox{in }\Omega,\\
   \mbox{div}\,u &=& 0\quad\mbox{in }\Omega,\label{eq:nse2}\\
   u&=& 0\quad\mbox{on }\partial\Omega.\label{eq:nse3}
\eeqal
We have the same space $X$ as for the Stokes problem in Section \ref{sec:stokes}.
The nonlinear operator reads
\beqa
   A(u,p) &=& \left(\begin{array}{c}(u\cdot\nabla)u+\frac12u\,\div u -\mu\Delta u-\nabla p \\
   			\div u\end{array}\right)
\eeqa
The  operator $A$ is the form (\ref{eq:spoanonlinear}) with $C(u):=(u\cdot\nabla)u+\frac12u\,\div u$, $Lu:=-\mu\Delta u$, $Bu:=\div u$
and $T:=0$.
In order to apply Proposition \ref{prop:noyetcoercive} we verify  conditions (A1), (A2),  (A3) and (A4):
\begin{description}
\item[\rm (A1)] $\langle Lu ,u\rangle_{V',V}=\mu\normV{u}^2$, and $L$ is continuous
$\langle Lu ,\phi\rangle_{V',V}\le \mu\normV{u}\normV{\phi}$.
\item[\rm (A2)] As in the Stokes case, this property follows from the Inf-sup condition for $B$.
\item[\rm (A3)] $\langle C(u) ,u\rangle_{V',V}=\langle (u\cdot\nabla)u+\frac12u\,\div u,u\rangle_{V',V}=0$.
\item[\rm (A4)] $\normVS{C(u)}\le c_N\normV{u}^2$ with a constant $c_N>0$ from Sobolev embedding.
\end{description}
Hence, with Proposition \ref{prop:noyetcoercive} we get for any $f\in H^{-1}(\Omega)$ a solution $u\in V$, $p\in Q$, of the stationary Navier-Stokes equations (\ref{eq:nse1})-(\ref{eq:nse3}).

\subsection{\bf Application to stabilized finite element approximations of the Navier-Stokes equations}

Now, we consider finite element discretization as in Section \ref{sec:stokesfem}  for the Navier-Stokes system.
In order to  address the lack of a finite dimensional analog of an inf-sup condition, we use a linear term $T\in L(Q_h,Q_h')$.
To  circumvent instabilities and numerical artifacts in the case of dominant convection, we allow for a further  
stabilization term $S:V_h\to V_h'$.
Then, the finite-dimensional approximation can be written in the following form: We seek $(u,p)\in X_h$ s.t. for any  $ (v,q)\in X_h$ it holds
\beqal\label{eq:nsediscrete}
     \langle (u\cdot\nabla)u+\tfrac12u\,\div u-\mu\Delta u+\nabla p+S(u),v\rangle_{V',V} &&\\
     + \langle \div u+Tp,q\rangle_{Q_h',Q_h}&=& \langle f,v\rangle_{V',V}.
\eeqal
As for the Stokes system, the operator $T$ can be chosen as any of the possibilities in Section \ref{sec:bps}, \ref{sec:lps}, or \ref{sec:bhs}.
For $S$ we also have several possibilities, e.g. by LPS:
\beqal\label{eq:lps2}&&
    \langle S(u),v\rangle_{V_h',V_h} :=
    \sum_{K\in \mathcal T_{2h}}\delta_K^{lps}(u) \int_K \left[\nabla u-\pi_h(\nabla u)\right]\!\cdot\! \left[\nabla v-\pi_h(\nabla v)\right],
\eeqal
with a function $\delta_K^{lps}:V_h\to\nR_{\ge0}$, depending on $u$ and $h_K$ on each patch $K\in\mathcal T_{2h}$, 
and an upper bound
\beqal\label{eq:upperboundlps}
   |\delta_K^{lps}(u)| &\le& c\| \nabla u\|_{L^2(K)}.
\eeqal
Alternatively, we may use the following stabilization of SUPG-type:
\beqal\label{eq:supg}
    \langle S(u),v\rangle_{V_h',V_h} &:= & 
    \sum_{K\in \mathcal T_{h}}\delta_K^{supg}(u) \int_K (u\cdot\nabla) u\cdot (u\cdot\nabla )v,
\eeqal
with $0\le \delta_K^{supg}(u)\le ch_K\|u\|_{L^\infty(K)}^{-1}$.

\begin{corollary}
Let $X_h$ the discrete space of equal-order elements of order $k$ described above. 
The Navier-Stokes equation (\ref{eq:nsediscrete})
with $f\in H^{-1}(\Omega)$,
$T\in L(Q_h,Q_h')$ from Section \ref{sec:bps}, \ref{sec:lps}, or \ref{sec:bhs}, in combination with  $S:V_h\to V_h'$
in (\ref{eq:lps2}) or (\ref{eq:supg}), has a  solution  $(u,p)\in X_h$.
\end{corollary}
\begin{proof}
In order to apply Proposition \ref{prop:noyetcoercive} we have to verify properties (A1), (A2a), (A2b), (A3) and (A4):
\begin{itemize}
\item (A1): The Laplace operator $L=-\mu\Delta$ is continuous and coercive.
\item (A2a): All pressure stabilizations $T\in L(Q_h,Q_h')$ in  Section \ref{sec:bps}, \ref{sec:lps}, and \ref{sec:bhs}
satisfy $ \langle Tp,p\rangle_{Q_h',Q_h}\ge 0$ for any $p\in Q_h$.
\item (A2b): This property was shown in Section \ref{sec:bps}, \ref{sec:lps}, and \ref{sec:bhs}.
\item (A3): For the nonlinear operator $C(u):=(u\cdot\nabla)u+\tfrac12u\,\div u-\mu\Delta u+S(u)$ it holds
\beqa
    \langle C(u),u\rangle_{V',V} &=&   \langle S(u),u\rangle_{V',V} \ge 0\qquad\forall u\in V_h.
\eeqa
\item (A4): According to Section \ref{sec:navierstokes} (the infinite-dimensional case), the Galerkin part $C(u)-S(u)$ satisfies (A4).
Hence, it is sufficient to verify $ \|S(u)\|_{V'} \le c_N \normV{u}^2$ for all $u\in V_h$. For the local projection (\ref{eq:lps2}), we use the $L^2$-stability of the
$L^2$-projection, $\|\pi_h(\nabla u)\|_{L^2(K)}\le c\| \nabla u\|_{L^2(K)}$, and (\ref{eq:upperboundlps}):
\beqa
 \langle S(u),v\rangle_{V',V} 
 &\le& \left(\sum_{K\in \mathcal T_{2h}}\delta_K^{lps}(u)^2 \| \nabla u-\pi_h(\nabla u)\|_{L^2(K)}^2\right)^{1/2}\\
 &&
   \cdot\|\nabla v-\pi_h(\nabla v)\|_{L^2(\Omega)}\\
   &\le& c\left(\sum_{K\in \mathcal T_{2h}}\delta_K^{lps}(u)^2 \| \nabla u\|_{L^2(K)}^2\right)^{1/2}\normV{v}\\
   &\le& c\left(\sum_{K\in \mathcal T_{2h}}\| \nabla u\|_{L^2(K)}^4\right)^{1/2}\normV{v}\
   \le\ c \normV{u}^2 \normV{v}.
\eeqa
This yields $\|S(u)\|_{V'} \le c_N \normV{u}^2$ for (\ref{eq:lps2}). The same bound we get for  (\ref{eq:supg}) with
H\"older's and Poincare's inequality:
\beqa
 \langle S(u),v\rangle_{V',V} 
 &\le&   \sum_{K\in \mathcal T_{h}}\delta_K^{supg}(u) \| (u\cdot\nabla) u\|_{L^1(K)} \|(u\cdot\nabla )v\|_{L^\infty(K)}\\
 &\le&   c\sum_{K\in \mathcal T_{h}} \| h_Ku\|_{L^2(K)}\|\nabla u\|_{L^2(K)} \|\nabla v\|_{L^\infty(K)}\\
 &\le&   c\left(\sum_{K\in \mathcal T_{h}} \| h_Ku\|^2_{L^2(K)} \right)^{1/2}\ \|\nabla u\|_{L^2(\Omega)}  \|\nabla v\|_{L^\infty(\Omega)}\\
 &\le&c \| hu\|_{L^2(\Omega)}\normV{u} \normV{v}.
\eeqa
Due to the finite dimension of $V_h$, we have norm equivalence $ \|\nabla v\|_{L^\infty(\Omega)}\le c \|\nabla v\|_{L^2(\Omega)}=c\normV{v}$
for all $v\in V_h$ (the constant $c$ is independent of $h$, because the arising norm contain exactly one derivative).
By an inverse estimate we have $\| hu\|_{L^2(\Omega)}\le c \normV{v}$. This yields $\normVS{S(u)}\le c\normV{u}^2$ and (A4).
\end{itemize}

\end{proof}

\normalsize

\bibliographystyle{plain} 
\bibliography{mybib}

\end{document}